%% file: ms.tex
\documentclass[
]{ifacconf}

\usepackage{graphicx}      
\usepackage[%
    square%
]{natbib}        
\input{preamble}
\begin{document}
\begin{frontmatter}

\title{Backstepping-based tracking control of the vertical gradient freeze crystal growth process}
\subtitle{\small This work has been submitted to IFAC for possible publication}

\thanks{%
    The work of the first author was partially funded by the Deutsche Forschungsgemeinschaft
    (DFG) under project number WI 4412/1-1.
}

\author[tud]{Stefan Ecklebe} 
\author[jku]{Nicole Gehring} 

\address[tud]{Institute of Control Theory, Technische Universität Dresden, Germany, (e-mail: stefan.ecklebe@tu-dresden.de)}
\address[jku]{Institute of Automatic Control and Control Systems Technology, Johannes Kepler University Linz, Austria,\\ (e-mail: nicole.gehring@jku.at)}

\begin{abstract}                
    The vertical gradient freeze crystal growth process is the main technique
    for the production of high quality compound semiconductors that are
    vital for today's electronic applications.
    A simplified model of this process consists of two 1D diffusion equations with free boundaries for the temperatures in crystal and melt.
    Both phases are coupled via an ordinary differential equation that describes the evolution of the moving solid/liquid interface.
    The control of the resulting two-phase Stefan problem is the focus of this contribution.
    A flatness-based feedforward design is combined with a multi-step backstepping approach to obtain a controller that tracks a reference trajectory for the position of the phase boundary.
    Specifically, based on some preliminary transformations to map the model into a time-variant PDE-ODE system, consecutive decoupling and backstepping transformations are shown to yield a stable closed loop.
    The tracking controller is validated in a simulation that considers the 
    actual growth of a Gallium arsenide single crystal.
\end{abstract}


\end{frontmatter}

\input{introduction} 
\input{modelling} 
\input{feedforward} 
\input{feedback} 
\input{simulation} 

\section{Concluding remarks}

Based on the very promising numerical results in Section \ref{sec:simresults} and even in the case of parametric uncertainties that was not shown, one would hope to test the controller in experiments.
Moreover, future works should include the detailed solution of the decoupling and kernel equations (see \eqref{eq:decoup_equations} and \eqref{eq:kernel_equations}) to clarify details omitted here for brevity.
In that context, the approach in \cite{Jadachowski2012,Freudenthaler2016} for the solution of time-varying kernel equations might help to improve the numerical solution of \eqref{eq:kernel_equations}.

\begin{ack}
The authors thank Simon Kerschbaum, formerly University of Erlangen–Nuremberg, for valuable discussions.
\end{ack}

\bibliography{bibliography}             

\end{document}

%% file: preamble.tex
\usepackage[utf8]{inputenc}
\usepackage[T1]{fontenc}
\usepackage[final]{microtype}

\usepackage{tudcolors}
\usepackage{graphicx}
\usepackage{tikz}
\usetikzlibrary{
	plotmarks,
	math,
	positioning,
	shapes,
	arrows,
	arrows.meta,
	backgrounds,
	circuits.logic.IEC,circuits.ee.IEC,
	patterns,
	decorations,
	decorations.pathmorphing,
	shapes.geometric,
	calc,
	fit,
	matrix
}
\graphicspath{%
	{include/}%
	{../include/}%
}

\usepackage{amsmath}
\usepackage{amssymb}
\usepackage[exponent-product=\cdot, tight-spacing=false]{siunitx}
\usepackage{mathtools}
\usepackage[fleqn]{nccmath}

\makeatletter
\DeclareRobustCommand\vdots{%
  \mathpalette\@vdots{}%
}
\newcommand*{\@vdots}[2]{%
  \sbox0{$#1\cdotp\cdotp\cdotp\m@th$}%
  \sbox2{$#1.\m@th$}%
  \vbox{%
    \dimen@=\wd0 %
    \advance\dimen@ -3\ht2 %
    \kern.5\dimen@
    \dimen@=\wd2 %
    \advance\dimen@ -\ht2 %
    \dimen2=\wd0 %
    \advance\dimen2 -\dimen@
    \vbox to \dimen2{%
      \offinterlineskip
      \copy2 \vfill\copy2 \vfill\copy2 %
    }%
  }%
}
\DeclareRobustCommand\ddots{%
  \mathinner{%
    \mathpalette\@ddots{}%
    \mkern\thinmuskip
  }%
}
\newcommand*{\@ddots}[2]{%
  \sbox0{$#1\cdotp\cdotp\cdotp\m@th$}%
  \sbox2{$#1.\m@th$}%
  \vbox{%
    \dimen@=\wd0 %
    \advance\dimen@ -3\ht2 %
    \kern.5\dimen@
    \dimen@=\wd2 %
    \advance\dimen@ -\ht2 %
    \dimen2=\wd0 %
    \advance\dimen2 -\dimen@
    \vbox to \dimen2{%
      \offinterlineskip
      \hbox{$#1\mathpunct{.}\m@th$}%
      \vfill
      \hbox{$#1\mathpunct{\kern\wd2}\mathpunct{.}\m@th$}%
      \vfill
      \hbox{$#1\mathpunct{\kern\wd2}\mathpunct{\kern\wd2}\mathpunct{.}\m@th$}%
    }%
  }%
}
\makeatother

\input{makros}

\usepackage{natbib}        

\usepackage{booktabs}

\usepackage{glossaries}
\newacronym{dps}{DPS}{distributed parameter system}
\newacronym{bvp}{BVP}{boundary value problem}
\newacronym{ivp}{IVP}{initial value problem}
\newacronym{fbp}{FBP}{free boundary problem}
\newacronym{sp}{SP}{Stefan problem}
\newacronym{opsp}{OPSP}{one-phase Stefan problem}
\newacronym{tpsp}{TPSP}{two-phase Stefan problem}
\newacronym{bc}{BC}{boundary condition}
\newacronym{sc}{SC}{Stefan condition}
\newacronym{ode}{ODE}{ordinary differential equation}
\newacronym{odes}{ODES}{ordinary differential equation system}
\newacronym{pde}{PDE}{partial differential equation}
\newacronym{fem}{FEM}{finite element method}
\newacronym{cfd}{CFD}{computational fluid dynamic}
\newacronym{vgf}{VGF}{Vertical Gradient Freeze}
\newacronym{tmf}{TMF}{travelling magnetic field}
\newacronym{gaas}{GaAs}{Gallium arsenide}
\newacronym{inp}{InP}{Indium phosphide}
\newacronym{dof}{DOF}{degree of freedom}
\newacronym{iss}{ISS}{Input State Stability}
\newacronym{pi}{CSI}{Poincaré inequality}
\newacronym{csi}{CSI}{Cauchy-Schwarz inequality}
\newacronym{yi}{YI}{Young's inequality}
\newacronym{lti}{LTI}{Linear Time Invartiant}
\newacronym{ltv}{LTV}{Linear Time Vartiant}
\newacronym{lqe}{LQE}{Linear Quadratic Estimator}
\newacronym{frde}{FDRE}{Filtering Riccati Differential Equation}

%% file: makros.tex
\DeclareMathOperator{\diag}{diag}

\newcommand{\tdiff}[1]{
	\tfrac{\mathrm{d}}{\mathrm{d}#1}
}

\newcommand{\dop}[1]{
	\,\mathrm{d}#1
}

\DeclarePairedDelimiterX{\norm}[1]{\lVert}{\rVert}{#1}

\newcommand{\bs}[1]{\boldsymbol{#1}}

\newcommand{\gen}{%
	\ensuremath{i}
}
\newcommand{\agen}{%
	\ensuremath{3-i}
}

\newcommand{\zs}{%
    \ensuremath{z}
}
\newcommand{\zl}{%
    \ensuremath{z}
}

\newcommand{\zb}{%
	\ensuremath{\Gamma}
}
\newcommand{\zbs}{%
	\ensuremath{\mathnormal{\zb}_{1}}
}
\newcommand{\zbl}{%
	\ensuremath{\mathnormal{\zb}_{2}}
}
\newcommand{\zbg}{%
	\ensuremath{\mathnormal{\zb}_{\gen}}
}

\newcommand{\errcomplete}[3]{%
	\ensuremath{\epsilon_{\mathrm{#1}}^{#2}%
		\ifthenelse{\equal{#3}{}}{}{(#3)}%
	}
}

\newcommand{\zisym}{%
	\ensuremath{\gamma}
}
\newcommand{\zicomplete}[4][]{%
	\ensuremath{%
			\zisym%
			_{\mathrm{#2}}^{#3}%
		\ifthenelse{\equal{#4}{}}{}{(#4)}%
	}
}
\newcommand{\vicomplete}[3]{%
	\ensuremath{\dot{\zisym}_{\mathrm{#1}}^{#2}%
		\ifthenelse{\equal{#3}{}}{}{(#3)}%
	}
}
\newcommand{\aicomplete}[3]{%
	\ensuremath{\ddot\zisym_{\mathrm{#1}}^{#2}%
		\ifthenelse{\equal{#3}{}}{}{(#3)}%
	}
}

\newcommand{\zi}[1][t]{%
	\ensuremath{\zicomplete{}{}{#1}}
}

\newcommand{\vi}[1][t]{%
	\ensuremath{\vicomplete{}{}{#1}}
}

\newcommand{\zir}{%
	\ensuremath{\zicomplete{r}{}{t}}
}
\newcommand{\vir}{%
	\ensuremath{\vicomplete{r}{}{t}}
}

\newcommand{\tm}{%
	\ensuremath{T_{\mathrm{m}}}%
}

\newcommand{\torig}{%
	\ensuremath{T}%
}

\newcommand{\torigz}{%
	\ensuremath{\torig(z, t)}%
}

\newcommand{\ts}{%
	\ensuremath{T_1}%
}
\newcommand{\tsz}{%
	\ensuremath{\ts(\zs, t)}%
}
\newcommand{\tl}{%
	\ensuremath{T_2}%
}
\newcommand{\tlz}{%
	\ensuremath{\tl(\zl, t)}%
}

\newcommand{\trefsz}{%
	\ensuremath{T_{1,\mathrm{r}}(\zs, t)}%
}
\newcommand{\treflz}{%
	\ensuremath{T_{2,\mathrm{r}}(\zs, t)}%
}

\newcommand{\zf}{%
    \ensuremath{z}
}

\newcommand{\tfixs}{%
    \ensuremath{\bar{\torig}_{1}}%
}
\newcommand{\tfixl}{%
    \ensuremath{\bar{\torig}_{2}}%
}
\newcommand{\tfixg}{%
    \ensuremath{\bar{\torig}_{\gen}}%
}

\newcommand{\tfixrefg}{%
    \ensuremath{\bar{\torig}_{\gen,\mathrm{r}}}%
}

\newcommand{\tfixgz}{%
	\ensuremath{\tfixg(\zf, t)}%
}

\newcommand{\tfixrefgz}{%
	\ensuremath{\tfixrefg(\zf, t)}%
}

\newcommand{\tferr}{%
	\ensuremath{\Delta\bar T}%
}
\newcommand{\tferrs}{%
	\ensuremath{\tferr_{1}}%
}
\newcommand{\tferrl}{%
	\ensuremath{\tferr_{2}}%
}
\newcommand{\tferrg}{%
	\ensuremath{\tferr_{\gen}}%
}
\newcommand{\tferrag}{%
	\ensuremath{\tferr_{\agen}}%
}

\newcommand{\tferrgz}{%
	\ensuremath{\tferrg(\zf, t)}%
}

\newcommand{\tfherr}{%
	\ensuremath{\Delta\check T}%
}
\newcommand{\tfherrs}{%
	\ensuremath{\tfherr_1}%
}
\newcommand{\tfherrl}{%
	\ensuremath{\tfherr_2}%
}
\newcommand{\tfherrg}{%
	\ensuremath{\tfherr_{\gen}}%
}
\newcommand{\tfherrag}{%
	\ensuremath{\tfherr_{\agen}}%
}

\newcommand{\tfherrsz}{%
	\ensuremath{\tfherrs(\zf, t)}%
}
\newcommand{\tfherrlz}{%
	\ensuremath{\tfherrl(\zf, t)}%
}
\newcommand{\tfherrgz}{%
	\ensuremath{\tfherrg(\zf, t)}%
}

\newcommand{\dzi}{%
	\ensuremath{\Delta\zi}%
}

\newcommand{\dvi}{%
	\ensuremath{\Delta\vi}%
}

\newcommand{\hcftrans}{%
	\mathnormal{\Phi}
}
\newcommand{\hcftransi}{%
	\hcftrans^{-1}
}
\newcommand{\hcftransg}{%
	\hcftrans_{\gen}
}

\newcommand{\hcftransig}{%
	\hcftransi_{\gen}
}

\newcommand{\hcftransgz}{%
	\hcftransg(\zf, t)
}

\newcommand{\hcftransigz}{%
	\hcftransig(\zf, t)
}

\newcommand{\hc}{%
	\ensuremath{\kappa}
}
\newcommand{\hcs}{%
	\ensuremath{\hc_1}
}
\newcommand{\hcl}{%
	\ensuremath{\hc_2}
}
\newcommand{\hcg}{%
	\ensuremath{\hc_{\gen}}
}

\newcommand{\hd}{%
	\ensuremath{\alpha}
}
\newcommand{\hds}{%
	\ensuremath{\hd_{1}}
}
\newcommand{\hdl}{%
	\ensuremath{\hd_{2}}
}
\newcommand{\hdg}{%
	\ensuremath{\hd_{\gen}}
}

\newcommand{\dn}{%
	\ensuremath{\rho}
}
\newcommand{\dmelt}{%
	\ensuremath{\dn_{\mathrm{m}}}
}

\newcommand{\bfd}{%
    \ensuremath{\left(-1\right)}
}
\newcommand{\bfdg}{%
    \bfd^{\gen}
}

\newcommand{\inpg}{%
	\ensuremath{u_{\gen}(t)}
}
\newcommand{\inps}{%
	\ensuremath{u_{1}(t)}
}
\newcommand{\inpl}{%
	\ensuremath{u_{2}(t)}
}
\newcommand{\inprs}{%
	\ensuremath{u_{1,\mathrm{r}}(t)}
}
\newcommand{\inprl}{%
	\ensuremath{u_{2,\mathrm{r}}(t)}
}
\newcommand{\inprg}{%
	\ensuremath{u_{\gen,\mathrm{r}}(t)}
}

\newcommand{\sode}{%
	\ensuremath{\bs{\xi}}
}
\newcommand{\spde}{%
	\ensuremath{\bs{x}}
}
\newcommand{\spded}{%
	\ensuremath{\bar{\bs{x}}}
}
\newcommand{\spdet}{%
	\ensuremath{\tilde{\bs{x}}}
}

\def\highlight<#1>{%
	\temporal<#1>{\color{HKS41}}{\color{HKS44}}{\color{HKS41}}
}

\newcommand{\scsym}{%
	\ensuremath{s}
}
\newcommand{\scfsym}{%
	\ensuremath{\bar{\scsym}}
}

\newcommand{\scfg}{%
	\ensuremath{\scfsym_{\gen}(\zi)}
}
\newcommand{\scfs}{%
	\ensuremath{\scfsym_{1}(\zi)}
}
\newcommand{\scfl}{%
	\ensuremath{\scfsym_{2}(\zi)}
}

\newcommand{\scg}{%
	\ensuremath{\scsym_{\gen}(t)}
}
\newcommand{\scs}{%
	\ensuremath{\scsym_{1}(t)}
}
\newcommand{\scl}{%
	\ensuremath{\scsym_{2}(t)}
}
\newcommand{\scj}{%
	\ensuremath{\scsym_j(t)}
}

\newcommand{\es}{%
	\ensuremath{q_\mathrm{p}}
}

\newcommand{\inpes}{%
	\ensuremath{\Delta u_{1}(t)}
}
\newcommand{\inpel}{%
	\ensuremath{\Delta u_{2}(t)}
}
\newcommand{\inpeg}{%
	\ensuremath{\Delta u_{\gen}(t)}
}

\newcommand{\inpev}{%
	\ensuremath{\bs{\Delta u}(t)}
}

%% file: introduction.tex
\section{Introduction}

\glsunset{pde}
\glsunset{ode}
From their invention in the early 1950s, the share of compound 
semiconductors such as \gls{gaas} in optoelectronic devices or
high speed digital applications has been constantly increasing.
This can be attributed to their larger bandgap and higher electron mobility 
in comparison to silicon (see e.g.~\cite{Jurisch2015}).
However, controlling the growth process is significantly more challenging than growing silicon crystals.
Single crystals of compound semiconductors are grown in crucibles via so-called Bridgman methods, of which the \gls{vgf}
method is a prominent member.
The technique works as follows (see \cite{Jurisch2015}):
For every charge, chunks of the material are molten in the crucible,
except for the so-called seed crystal at the bottom, which is used to establish
the desired growth orientation.
Then, heaters inside the furnace are used to create a vertical temperature gradient in the material, that, by moving the temperature profile from bottom to top, lets the crystal grow accordingly.
Thus, the position of the interface between the solid and the still liquid phase evolves in time, resulting in a \gls{sc}.
This moving interface 
causes the domains for crystal and melt to change over time.
That renders the process a free boundary problem, which is 
inherently nonlinear and known as the \gls{tpsp}, see e.g.~\cite{Crank84}.

The \gls{tpsp} is the main challenge in controlling the \gls{vgf} process, which is why most research focuses solely on this (abstract) aspect.
An overview of control strategies for Stefan problems can be found in \cite{Koga2022}.
Specifically, methods include the energy-based designs in \cite{petrus2010,Koga2020}, 
a flatness-based approach in \cite{Ecklebe2021} and the optimal control designs in \cite{Kang1995, Hinze2009,Baran2018,Ecklebe2021}.
Considering only a \gls{opsp}, where the temperature in one phase is assumed to be at melting temperature everywhere,  \cite{Koga2019a} uses a backstepping design (see \cite{book:KrsticBackstepping}, in general) to stabilise a desired fixed position of the phase boundary.
However, ignoring the second phase only allows to correct the interface position in one direction, which is insufficient for the tracking control of the \gls{vgf} process.
Therefore, \cite{Ecklebe2020} suggests a backstepping design for the \gls{tpsp} that neglects the coupling between both phases and stabilises them independently of one another.
In the last years, new results on backstepping control have emerged that now allow to incorporate the couplings in the design.
Therefore, here, a tracking controller is designed for the \gls{tpsp}, with the feedforward design taken from \cite{Ecklebe2020}.
The suggested approach applies the multi-step, transformation-based design in \cite{Gehring2021} for so-called PDE-ODE systems.
In its final design step, a Volterra integral transformation based on the results for coupled parabolic \glspl{pde} in \cite{Kerschbaum2020} is used.
In contrast to \cite{Kerschbaum2020} and \cite{Gehring2021}, here, additional challenges arise due to the time-varying nature of the PDE-ODE system based on the \gls{tpsp}.

The paper is structured as follows:
Section \ref{sec:modelling} gives the mathematical model for the \gls{vgf} plant, i.e.\ the \gls{tpsp}.
Using a flatness-based approach, the feedforward part of the tracking control design is addressed in Section \ref{sec:feedforward}.
Linearisation of the \gls{vgf} model around this reference results in a time-varying \gls{pde}-\gls{ode} system.
For that, in Section \ref{sec:feedback}, the tracking controller is derived by means of a multi-step backstepping design.
The numerical results in Section \ref{sec:simulation} give insight into the closed-loop performance.

%% file: modelling.tex
\section{Modelling and problem statement}
\label{sec:modelling}

\begin{figure}
	\centering
	\def\svgwidth{0.9\linewidth}
	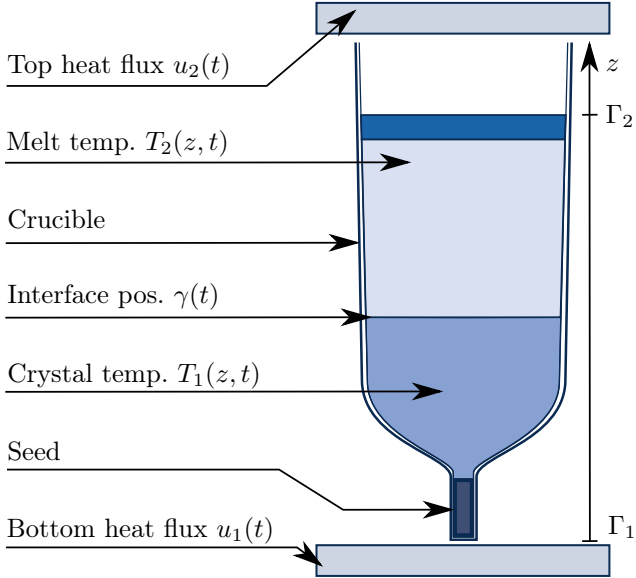
	\caption{
	    Cross section of a typical \gls{vgf} plant (without jacket heaters).
	    }
	\label{fig:vgf_plant}
\end{figure}


This section briefly revisits a simplified 1D model of a \gls{vgf} plant (see \cite{Ecklebe2021} for the comprehensive form) and specifies the control objective.

\subsection{Plant Model}
\label{ssc:plant_model}

We only focus on the processes inside the crucible, 
a cross section of which is depicted in Figure \ref{fig:vgf_plant}.
Furthermore, deviating from \cite{Ecklebe2021} the jacket heaters are 
assumed to act as active insulation and therefore ignored.
In the crucible, the main dynamics are given by the heat diffusion processes
in the solid crystal and the liquid melt since the convection in the melt is negligible due to the low Prandtl number of the
considered materials \cite[Tab.~9.1]{Jurisch2015}.
In general, the heat transport processes have to be modelled in 3D.
However, 
by exploiting the rotational symmetry of the crucible,
the domain can be reduced to two dimensions.
Furthermore, by assuming an overall slow growth regime,
the radial temperature gradients can be neglected,
which yields a 1D temperature profile $\torigz$ in the spatial coordinate $z$
and temporal coordinate $t$.
Yet, a closed-form description of the diffusion process 
in the complete system
is not feasible, as the physical properties of the material change abruptly
with the phase transition between liquid and solid.
Therefore, the overall temperature profile $\torigz$ is split into two parts,
$\tsz$ in the crystal and $\tlz$ in the melt,
both governed by the diffusion equations \eqref{eq:lin_heat_eq_s_pde} and \eqref{eq:lin_heat_eq_l_pde} in 
\begin{subequations}
	\label{eq:lin_heat_eqs}
	\begin{align}
		\label{eq:lin_heat_eq_s_pde}
		\partial_t \tsz &= \hds\partial_{\zs}^2 \tsz,
		    \quad \zs \in (\zbs, \zi)
		\\
		\label{eq:lin_heat_eq_s_bc1}
		\hcs\partial_{\zs} \ts (\zbs, t) &= -\inps\\
		\label{eq:lin_heat_eq_s_bc2}
		\ts(\zi, t) &= \tm \\
		\label{eq:lin_heat_eq_l_pde}
		\partial_t \tlz &= \hdl\partial_{\zl}^2 \tlz,
		    \quad \zl \in (\zi, \zbl)
		\\
		\label{eq:lin_heat_eq_l_bc1}
		\hcl\partial_{\zl} \tl (\zbl, t) &= \inpl\\
		\label{eq:lin_heat_eq_l_bc2}
		\tl(\zi, t) &= \tm,
	\end{align}
    with constant parameters for the 
    heat diffusivity $\hdg$  and heat conductivity $\hcg$ where $i=1,2$.
    Therein, due to the ongoing phase transition 
    at the moving interface $\zi$, the
    Dirichlet \glspl{bc} \eqref{eq:lin_heat_eq_s_bc2} and \eqref{eq:lin_heat_eq_l_bc2} correspond to the melting temperature $\tm$.
    The system is actuated via the heat flows originating from the bottom and top heaters
    $\inps$ and $\inpl$ at the lower and upper system boundaries $\zbs$ and $\zbl$, 
    respectively.
    This results in the
    Neumann \glspl{bc} \eqref{eq:lin_heat_eq_s_bc1} and \eqref{eq:lin_heat_eq_l_bc1}.
    The evolution of the interface itself is driven by the heat flows from both phases
    into the interface as well as the latent heat released during the growth process.
    Denoting the density at melting temperature by $\dmelt$ and the specific heat of solidification by $\es$, an energy balance yields the Stefan condition
    \begin{equation}
        	\label{eq:stefan_cond_sep}
        	\dmelt\es \vi = \hcs \partial_z \ts(\zi, t) 
    	    - \hcl \partial_z \tl(\zi, t).
    \end{equation}
\end{subequations} 
The system \eqref{eq:lin_heat_eqs} is defined for $t\ge0$, with appropriate initial conditions for $\tsz$, $\tlz$ as well as $\zi$.
It describes a \gls{tpsp}, which is a nonlinear free boundary initial value problem, as the solution domains depend on $\zi$.

\subsection{Problem Statement}
\label{ssc:problem_statement}

From a technological point of view, the system \eqref{eq:lin_heat_eqs} has to be controlled such that a usable crystal is grown.
For that, in particular, the position $\zi$ of the phase boundary should follow a prescribed reference $\zir$, thus growing the crystal to a desired length.
Meanwhile, it is of vital importance to control the temperature gradient in the crystal, especially the gradient $\partial_{\zf}\ts(\zi, t)$ at the interface, in order to avoid cracking and to achieve a low dislocation density (and thus a high quality) of the resulting
crystal (see \cite{Vanhellemont2013}).
In what follows, the backstepping technique is used to stabilise the system \eqref{eq:lin_heat_eqs} along a desired reference trajectory by means of the control inputs $\inps$ and $\inpl$, i.e.\ the bottom and top heaters.
The following assumption is imposed for the control design:

\begin{assum}
\label{ass:interface}
    The position of the interface satisfies $\zi \in (\zbs, \zbl)$, $t\ge0$.
\end{assum}

This assumption is technologically motivated. It is always met since the initial seeding guarantees $\zi>\zbs$ for the position of the interface and not all melt is used for the crystallisation, resulting in $\zi<\zbl$.

%% file: include/kronos_overview_text_asym_mirrored_2.pdf_tex
\begingroup%
  \makeatletter%
  \providecommand\color[2][]{%
    \errmessage{(Inkscape) Color is used for the text in Inkscape, but the package 'color.sty' is not loaded}%
    \renewcommand\color[2][]{}%
  }%
  \providecommand\transparent[1]{%
    \errmessage{(Inkscape) Transparency is used (non-zero) for the text in Inkscape, but the package 'transparent.sty' is not loaded}%
    \renewcommand\transparent[1]{}%
  }%
  \providecommand\rotatebox[2]{#2}%
  \newcommand*\fsize{\dimexpr\f@size pt\relax}%
  \newcommand*\lineheight[1]{\fontsize{\fsize}{#1\fsize}\selectfont}%
  \ifx\svgwidth\undefined%
    \setlength{\unitlength}{253.64781189bp}%
    \ifx\svgscale\undefined%
      \relax%
    \else%
      \setlength{\unitlength}{\unitlength * \real{\svgscale}}%
    \fi%
  \else%
    \setlength{\unitlength}{\svgwidth}%
  \fi%
  \global\let\svgwidth\undefined%
  \global\let\svgscale\undefined%
  \makeatother%
  \begin{picture}(1,0.95722593)%
    \lineheight{1}%
    \setlength\tabcolsep{0pt}%
    \put(0,0){\includegraphics[width=\unitlength,page=1]{kronos_overview_text_asym_mirrored_2.pdf}}%
    \put(0.00256349,0.71127089){\color[rgb]{0,0,0}\makebox(0,0)[lt]{\lineheight{0}\smash{\begin{tabular}[t]{l}Melt temp. $T_{\mathrm{2}}(z,t)$\end{tabular}}}}%
    \put(0.00384089,0.32452272){\color[rgb]{0,0,0}\makebox(0,0)[lt]{\lineheight{0}\smash{\begin{tabular}[t]{l}Crystal temp. $T_{\mathrm{1}}(z,t)$\end{tabular}}}}%
    \put(0.00320175,0.83926558){\color[rgb]{0,0,0}\makebox(0,0)[lt]{\lineheight{0}\smash{\begin{tabular}[t]{l}Top heat flux $u_{\mathrm{2}}(t)$\end{tabular}}}}%
    \put(0.00315147,0.06722185){\color[rgb]{0,0,0}\makebox(0,0)[lt]{\lineheight{0}\smash{\begin{tabular}[t]{l}Bottom heat flux $u_{\mathrm{1}}(t)$\end{tabular}}}}%
    \put(0.00384089,0.19498603){\color[rgb]{0,0,0}\makebox(0,0)[lt]{\lineheight{0}\smash{\begin{tabular}[t]{l}Seed\end{tabular}}}}%
    \put(0.00384089,0.58037187){\color[rgb]{0,0,0}\makebox(0,0)[lt]{\lineheight{0}\smash{\begin{tabular}[t]{l}Crucible\end{tabular}}}}%
    \put(0.00384089,0.4556694){\color[rgb]{0,0,0}\makebox(0,0)[lt]{\lineheight{0}\smash{\begin{tabular}[t]{l}Interface pos. $\gamma(t)$\end{tabular}}}}%
    \put(0.99220468,0.84273818){\color[rgb]{0,0,0}\makebox(0,0)[lt]{\lineheight{0}\smash{\begin{tabular}[t]{l}$z$\end{tabular}}}}%
    \put(0.99626585,0.07143877){\color[rgb]{0,0,0}\makebox(0,0)[lt]{\lineheight{0}\smash{\begin{tabular}[t]{l}$\Gamma_{\mathrm{1}}$\end{tabular}}}}%
    \put(0.99220468,0.75595291){\color[rgb]{0,0,0}\makebox(0,0)[lt]{\lineheight{0}\smash{\begin{tabular}[t]{l}$\Gamma_{\mathrm{2}}$\end{tabular}}}}%
  \end{picture}%
\endgroup%

%% file: feedforward.tex
\section{Flatness-based feedforward design}
\label{sec:feedforward}

The flatness-based feedforward design for \eqref{eq:lin_heat_eqs} can be found in \cite{Ecklebe2020} and is briefly revisited here to make the paper self-contained. The design is based on the general results in \cite{DPRM03,Rudolph2003} and makes use of a power series representation of the system solutions in terms of a so-called flat output.
For \eqref{eq:lin_heat_eqs}, the flat output  
$\bs{y}(t) \coloneqq [\zi,\partial_{\zf}\ts(\zi, t)]^{\top}$ comprises the interface position as well as the temperature gradient at the interface, i.e.\ the key quantities to be controlled (see Section~\ref{ssc:problem_statement}).
Based on the flatness-based parametrisation of the solutions, choosing a desired reference trajectory $\bs{y}_\text{r}(t) \coloneqq [\zir,\partial_{\zf}T_\text{1,r}(\zi, t)]^{\top}$ for the flat output enables the computation of the corresponding temperature profiles, which in turn imply a reference for the inputs by \eqref{eq:lin_heat_eq_s_bc1} and \eqref{eq:lin_heat_eq_l_bc1}.
Thus, the feedforward design yields the references
\begin{equation}
\label{eq:reference}
    \zir,\; 
    \trefsz,\; \treflz,\;
    \inprs,\; \inprl.
\end{equation}
Note that the reference $\bs{y}_\text{r}(t)$ is usually chosen to transition the system \eqref{eq:lin_heat_eqs} between two steady states (see also Section \ref{sec:simulation}). In view of the infinitely many time derivatives of $\bs{y}(t)$ occurring in the flatness-based parametrisation, the functions $\bs{y}_\text{r}(t)$ have to be of a certain Gevrey order, for the series to converge. 
As a result, all functions in \eqref{eq:reference} are of the same Gevrey order.

%% file: feedback.tex
\section{Backstepping-based feedback design} 
\label{sec:feedback}

In order to apply the systematic, multi-step backstepping design in \cite{Gehring2021} to the \gls{vgf} plant, first, \eqref{eq:lin_heat_eqs} is mapped into a form consistent with \cite{Gehring2021}.
Mainly, this involves a transformation of $z$ to fix the spatial domain as well as a linearisation around the reference \eqref{eq:reference}.
The resulting \gls{pde}-\gls{ode} model has time-varying coefficients due to their dependence on the reference trajectories.
The control design is then split into two steps.
A decoupling transformation is used to stabilise the \gls{ode} and to decouple it from the \gls{pde}.
Using a backstepping transformation, the second step allows to stabilise the \gls{pde} by choice of the control input and results in a stable closed loop.

\subsection{Preliminary transformations}
\label{sec:prelim_trans}

Similar to \cite{Ecklebe2021}, the mappings
\begin{equation}
\label{eq:map_fixed}
    z \mapsto \sigma_i(z,\zi) = \frac{\zs - \zi}{\Gamma_i - \zi}, \qquad
    i=1,2
\end{equation}
are introduced in order to obtain time-invariant spatial domains with $z\in[0,1]$ for \eqref{eq:lin_heat_eqs}. 
Recall that \eqref{eq:map_fixed} is well-defined and invertible due to Assumption \ref{ass:interface}.
Thus, defining pointwise correspondence in 
$\bar{T}_i(z,t) = T_i(\sigma_i^{-1}(z, \zi),t)$, $i=1,2$, 
the model \eqref{eq:lin_heat_eqs} takes the form
\begin{subequations}
\label{eq:lin_heat_sys_fixed}
    \begin{align}
        \begin{split}
            \partial_t \tfixgz &= 
            \bar \lambda_{\gen}(\zi) \partial_{\zf}^2 \tfixgz\\
            &\qquad+ \bar \psi_{\gen}(\zf, \zi, \vi) \partial_{\zf} \tfixgz,
            \; \zf \in (0,1)
        \end{split}
        \\
        \tfixg(0, t) &= \tm \\
		\partial_{\zf}\tfixg(1, t) &= \bar q_{\gen}\inpg \\
        \vi &= \scfs\partial_{\zf}\tfixs(0, t) + \scfl\partial_{\zf}\tfixl(0, t), 
    \end{align}%
\end{subequations}%
with $\bar \lambda_{\gen}(\zi) \coloneqq \frac{\hdg}{\left(\zbg - \zi\right)^2}$,
$\bar \psi_{\gen}(\zf, \zi, \vi) \coloneqq \frac{(1 - \zf)\vi}{\left(\zbg - \zi\right)}$,
$\bar q_{\gen}(\zi) \coloneqq  \frac{\bfdg(\zbg - \zi)}{\hc}$ and
$\scfg \coloneqq -\frac{\bfdg\hcg}{ \dmelt\es(\zbg - \zi)}$
for $i=1,2$.

As the model \eqref{eq:lin_heat_sys_fixed} with the fixed spatial domain is nonlinear in the position $\zi$ of the interface, it is linearised around the reference trajectory \eqref{eq:reference}.
Introducing the errors 
$\tferrgz \coloneqq \tfixgz - \tfixrefgz$, $\dzi \coloneqq \zi - \zir$ and
$\inpeg \coloneqq \inpg - \inprg$
for $i=1,2$, the time-variant system%
\begin{subequations}
    \label{eq:sys_1d_fix_lin}
    \begin{align}
        \label{eq:sys_1d_fix_lin_pde}
        \begin{split}
            \partial_t \tferrgz &= 
                \lambda_{\gen}(t) \partial_{\zf}^2 \tferrgz
                + \psi_{\gen}(\zf, t) \partial_{\zf} \tferrgz
                \\&\;\quad
                + c_{\gen,\gen}(\zf, t)\partial_{\zf}\tferrg(0, t) 
                + r_{\gen}(\zf ,t)\dzi
               \\&\;\quad
                + c_{\gen,\agen}(\zf, t)\partial_{\zf}\tferrag(0, t),
                \; \zf \in (0,1)
        \end{split}\\
        \tferrg(0, t) &= 0 \\
        \partial_{\zf} \tferrg(1, t) &=  
            p_{\gen}(t)\dzi 
            + q_{\gen}(t)\inpeg 
            \\
        \begin{split}
            \dvi &= 
                \scs \partial_{\zf}\tferrs(0, t) 
                +\scl \partial_{\zf}\tferrl(0, t) 
               \\&\;\quad
                + d(t) \dzi
        \end{split}
    \end{align}
\end{subequations}
is obtained.
Therein,
\begin{align*}
	f_{\gen}(\zf ,t) &\coloneqq \textstyle\frac{2\hdg}{\left(\zbg - \zir\right)^3} \partial_{\zf}^2\tfixrefgz 
				+ \frac{(1 - \zf)\vir}{\left(\zbg - \zir\right)^2} \partial_{\zf}\tfixrefgz\\
	g_{\gen}(\zf, t) &\coloneqq \textstyle\frac{1 - \zf}{\zbg - \zir} \partial_{\zf} \tfixrefgz\\
	d(t) &\coloneqq \textstyle\sum_{i=1}^{2} \frac{\scg}{\zbg - \zir}\partial_{\zf}\tfixrefg(0, t)
\end{align*}
and
$\lambda_{\gen}(t) \coloneqq \bar \lambda_{\gen}(\zir)$,
$\psi_{\gen}(\zf, t)	\coloneqq \bar \psi_{\gen}(\zf, \zir, \vir)$,
$c_{i,j}(\zf ,t) \coloneqq g_{\gen}(\zf, t)\scj$, $j=1,2$,
$r_{\gen}(\zf ,t) \coloneqq f_{\gen}(\zf, t) + g_{\gen}(\zf, t) d(t)$,
$q_{\gen}(t) \coloneqq \bar q_{\gen}(\zir)$,
$p_{\gen}(t) \coloneqq -\frac{\bfdg\inprg}{\hcg}$ and
$\scg \coloneqq \scfsym_{\gen}(\zir)$ for $i=1,2$.

Finally, to simplify the backstepping controller design in the subsequent sections, the convective terms in \eqref{eq:sys_1d_fix_lin_pde} are eliminated.
Using the Hopf-Cole-type transformation
$\tferrgz = \hcftransgz\tfherrgz$
with
\begin{equation}
    \hcftransgz \coloneqq \exp \left(-\textstyle\int_0^{\zf}\frac{\psi_{\gen}(\zeta,t)}{2\lambda_{\gen}(t)} \dop \zeta\right),
\end{equation}
the system \eqref{eq:sys_1d_fix_lin} reads
\begin{subequations}
    \label{eq:sys_1d_fix_lin_hc}
    \begin{fleqn}
        \begin{align}
            \begin{split}
                \partial_t \tfherrgz &= 
                    \lambda_{\gen}(t) \partial_{\zf}^2 \tfherrgz
                    + \check a_{\gen}(\zf, t) \tfherrgz
                    \\&\;\quad
                    + \check c_{\gen,\gen}(\zf, t)\partial_{\zf}\tfherrg(0, t) 
                    + \check r_{\gen}(\zf ,t)\dzi 
                    \\&\;\quad
                    + \check c_{\gen,\agen}(\zf, t)\partial_{\zf}\tfherrag(0, t),
                    \; \zf \in (0,1)
            \end{split}\\
            \tfherrg(0, t) &= 0 \\
            \partial_{\zf} \tfherrg(1, t) &= 	
                \check{b}_{\gen}(t) \tfherrg(1, t)
                \!+ \check{p}_{\gen}(t) \dzi
                \!+ \check{q}_{\gen}(t) \inpeg
                \raisetag{.8\baselineskip}
                \\
            \begin{split}
                \dvi &= 
                    \scs \partial_{\zf}\tfherrs(0, t) 
                    +\scl \partial_{\zf}\tfherrl(0, t) 
                    \\&\;\quad
                    + d(t) \dzi
            \end{split}
        \end{align}
    \end{fleqn}
\end{subequations} 
for $i=1,2$, with
\begin{align*}	
	\check a_{\gen}(\zf, t)	&\coloneqq \hcftransigz \big(
		\lambda_{\gen}(t) \partial_{\zf}^2 \hcftransgz 
		+ \psi_{\gen}(\zf, t) \partial_{\zf} \hcftransgz 
		\nonumber\\&\hphantom{\coloneqq\hcftransigz\big(}
		- \partial_{t} \hcftransgz 
	\big)\\
	\check c_{i,j}(\zf,t) &\coloneqq \hcftransigz c_{i,j}(\zf,t), \qquad j=1,2\\
	\check r_{\gen}(\zf,t) &\coloneqq \hcftransigz r_{\gen}(\zf,t),\;
	\check b_{\gen}(t) \coloneqq 
	    -\hcftransig(1,t)\partial_{\zf}\hcftransg(1,t)\\ 
	\check q_{\gen}(t) &\coloneqq \hcftransig(1,t) q_{\gen}(t),\;
	\check p_{\gen}(t) \coloneqq \hcftransig(1,t) p_{\gen}(t).
\end{align*}

\subsection{Extended Plant}

The importance of tracking a reference $\partial_{\zf}T_{1,\text{r}}(\zi, t)$ for the gradient at the interface (cf.\ Section~\ref{ssc:problem_statement}) is incorporated in the control design by using a dynamic state feedback controller with the integral dynamics
\begin{equation}
\label{eq:dyn_ext}
    \dot\epsilon(t) = \partial_{\zf}\tfherrs(0, t).
\end{equation}
Thus, augmenting the system \eqref{eq:sys_1d_fix_lin_hc} by the dynamics \eqref{eq:dyn_ext} results in the extended system \begin{subequations}
    \label{eq:sys_1d_fix_lin_hc_vec}
    \begin{align}
        \label{eq:sys_1d_fix_lin_hc_vec_ode}
        \dot\sode(t) &= 
            \bs{F}(t) \sode(t) 
            + \bs{S}(t) \partial_{\zf}\spde(0,t)
            \\
        \label{eq:sys_1d_fix_lin_hc_vec_bc0}
        \spde(0, t) &= \bs{0}
            \\
        \label{eq:sys_1d_fix_lin_hc_vec_pde}
        \begin{split}
            \partial_t \spde(\zf,t) &=
                \bs{\Lambda}(t)\partial_{\zf}^2 \spde(\zf,t)
                + \bs{A}(\zf, t) \spde(\zf,t)
                \\&\qquad
                + \bs{C}(\zf, t)\partial_{\zf} \spde(0, t)
                + \bs{R}(\zf, t) \sode(t) 
        \end{split} \\ 
        \label{eq:sys_1d_fix_lin_hc_vec_bc1}
        \partial_{\zf} \spde(1, t) &= 
            \bs{B}(t) \spde(1,t)
            + \bs{P}(t) \sode(t)
            + \bs{Q}(t) \inpev,
    \end{align}
\end{subequations}
with the \gls{ode} state $\sode(t) \coloneqq [\dzi,\epsilon(t)]^{\top}$, the \gls{pde} state
$\spde(\zf,t) \coloneqq [\tfherrsz,\tfherrlz]^{\top}$
and the input $\inpev \coloneqq [\inpes,\inpel]^{\top}$ as well as matrices $\bs{F}(t) \coloneqq \diag(d(t), 0)$
$\bs{\Lambda}(t) \coloneqq \diag(\lambda_{1}(t), \lambda_{2}(t))$,
$\bs{A}(\zf, t) \coloneqq \diag(
    \check a_{1}(\zf, t), \check a_{2}(\zf,t))$,
$\bs{B}(t) \coloneqq \diag(
		\check b_{1}(t), \check b_{2}(t))$,
$\bs{Q}(t) \coloneqq \diag(\check q_{1}(t), \check q_{2}(t))$,
\begin{align*}
	\bs{S}(t) &\coloneqq \begin{bmatrix}
		\scs	& \scl \\
		1		& 0 \\
	\end{bmatrix},\;
	&\bs{C}(\zf, t) &\coloneqq \begin{bmatrix}
		c_{1,1}(\zf,t)	& c_{1,2}(\zf,t) \\
	    c_{2,1}(\zf,t)	& c_{2,2}(\zf,t) \\
	\end{bmatrix}
    \\
	\bs{R}(\zf, t) &\coloneqq \begin{bmatrix}
		\check r_{1}(\zf, t)	& 0 \\
		\check r_{2}(\zf, t)	& 0 \\
	\end{bmatrix},\;
	&\bs{P}(t) &\coloneqq \begin{bmatrix}
		\check p_{1}(t)		&0 \\
		\check p_{2}(t)		&0\\
	\end{bmatrix}.
\end{align*}

\begin{rem}
\label{rem:diffsort}
    In view of the definition of the diffusion coefficients $\lambda_1(t)$ and $\lambda_2(t)$ as functions of the reference $\zir$, they do not satisfy a specific ordering for all times $t$ in general, as is usually assumed for the backstepping design.
\end{rem}

\subsection{Decoupling transformations}
\label{ssc:decoupling}

Noting that the boundary value $\partial_{\zf}\spde(0,t)$ takes the role of an input w.r.t.\ the \gls{ode}  \eqref{eq:sys_1d_fix_lin_hc_vec_ode}, one would like to use a feedback $\partial_{\zf}\spde(0,t)=-\bs{K}(t)\sode(t)$ in order to render $\dot\sode(t) = \big(\bs{F}(t) - \bs{S}(t)\bs{K}(t)\big)\sode(t)$ stable.
While this feedback cannot be implemented as $\inpev$ is the control input, this objective motivates the decoupling transformation
\begin{equation}
    \label{eq:decoupl_trafo}
    \spdet(\zf,t) \coloneqq \spde(\zf,t) - \bs{N}(\zf, t)\bs{x}(t), \quad \bs{N}(\zf,t)\in\mathbb{R}^{2\times2}.
\end{equation}
It introduces error variables, such that $\spded(\zf,t)\to\bs{0}$ for $t\to\infty$ yields the desired virtual feedback for the \gls{ode} if $\partial_{\zf} \bs{N}(0, t) = -\bs{K}(t)$. By applying \eqref{eq:decoupl_trafo}, \eqref{eq:sys_1d_fix_lin_hc_vec} is mapped into
\begin{subequations}
    \label{eq:sys_1d_fix_dec}
    \begin{align}
        \label{eq:sys_1d_fix_dec_ode}
        \dot\sode(t) &= 
            \bar{\bs{F}} \sode(t) 
            + \bs{S}(t) \partial_{\zf}\spded(0,t)
            \\
        \label{eq:sys_1d_fix_dec_bc0}
        \spded(0, t) &= \bs{0}
            \\
        \label{eq:sys_1d_fix_dec_pde}
        \begin{split}
            \partial_t \spded(\zf,t) &=
                \bs{\Lambda}(t) \partial_{\zf}^2 \spded(\zf,t)
                + \bs{A}(\zf, t) \spded(\zf,t)
               \\&\qquad
                + \bar{\bs{C}}(\zf, t) \partial_{\zf} \spded(0, t)
        \end{split} \\
        \label{eq:sys_1d_fix_dec_bc1}
        \partial_{\zf} \spded(1, t) &= 
            \bs{B}(t) \spded(1,t)
            + \bar{\bs{P}}(t) \sode(t)
            + \bs{Q}(t) \inpev,
    \end{align}
\end{subequations}
where
\begin{subequations}
        \begin{align}
            \label{eq:def_barF}
        	\bar{\bs{F}} &\coloneqq \bs{F}(t) + \bs{S}(t)\bs{N}(0,t)\\
        	\bar{\bs{C}}(\zf, t) &\coloneqq \bs{C}(\zf, t) 
        		- \bs{N}(\zf, t) \bs{S}(t)
        		\\
        	\bar{\bs{P}}(t) &\coloneqq \bs{P}(t) 
        		+ \bs{B}(t)\bs{N}(1,t) 
        		- \partial_{\zf} \bs{N}(1,t)
        \end{align}
\end{subequations}
if $\bs{N}(\zf,t)$ satisfies the decoupling equations
\begin{subequations}
    \label{eq:decoup_equations}
    \begin{align}
        \label{eq:decoup_equations_pde}
        \partial_t \bs{N}(\zf, t) &= 
        \bs{\Lambda}(t) \partial_{\zf}^2 \bs{N}(\zf, t)
        + \bs{A}(\zf, t) \bs{N}(\zf, t)
        \nonumber\\&\quad
        - \bs{N}(\zf, t) \bar{\bs{F}}
        - \bs{C}(\zf, t) \bs{K}(t) + \bs{R}(\zf, t)
        \\
        \label{eq:decoup_equations_0}
        \bs{N}(0, t) &= \bs{0} 
        \\
        \label{eq:decoup_equations_0z}
        \partial_{\zf} \bs{N}(0, t) &= -\bs{K}(t).
    \end{align} 
\end{subequations}
Note that the conditions \eqref{eq:decoup_equations_pde} and \eqref{eq:decoup_equations_0} on $\bs{N}(\zf,t)$ follow from the objective to decouple the \gls{pde} \eqref{eq:sys_1d_fix_dec_bc0}--\eqref{eq:sys_1d_fix_dec_bc1} from the \gls{ode}, meaning that the \gls{ode} state $\sode(t)$ only impacts the actuated boundary \eqref{eq:sys_1d_fix_dec_bc1}, where it could be compensated by choice of $\inpev$ in view of the invertibility of $\bs{Q}(t)$.
This property of $\bs{Q}(t)$ also allows to choose the design parameters in $\bs{K}(t)$ such that $\bar{\bs{F}}$ is a time-invariant Hurwitz matrix.
Consequently, \eqref{eq:sys_1d_fix_dec_ode} is stabilised.
In particular, $\sode(t)\to\bs{0}$ for $t\to\infty$ if $\partial_{\zf}\spdet(0,t)\to\bs{0}$. The latter is guaranteed by the backstepping step in Section \ref{ssc:backstepping}.

For brevity, the solution of the \gls{ivp} \eqref{eq:decoup_equations} is only sketched here.
Basically, using the eigenvectors of a constant matrix $\bar{\bs F}$ allows to rewrite the matrix-valued \gls{ivp} as two projected, vector-valued ones.
Then, a power series ansatz is used, similar to the derivation of a flatness-based parametrisation for \eqref{eq:lin_heat_eqs}.
In the end, the solution of \eqref{eq:decoup_equations} can be written in the form $\bs{N}(\zf,t)=\sum_{j=0}^{\infty}\bs{L}_j(t)\frac{\zf^j}{j!}$, where the matrices $\bs{L}_j(t)\in\mathbb{R}^{2\times2}$ are recursively defined based on $\bs{K}(t)$, $\bs{\Lambda}(t)$, $\bs{A}(\zf,t)$, $\bs{C}(\zf,t)$ and $\bs{R}(\zf,t)$ as well as their time derivatives.
Noting that the entries of all of these matrices are of a certain Gevrey order due to their dependence on the reference \eqref{eq:reference}, an appropriate choice of $\bs{K}(t)$ ensures the convergence of the series and the existence of a solution $\bs{N}(\zf,t)$.

\subsection{Backstepping transformation}
\label{ssc:backstepping}

Following the stabilization of the \gls{ode}, the Volterra integral transformation
\begin{equation}
    \label{eq:backstepp_trafo}
    \spdet(\zf,t) 
        \!\coloneqq\! \mathcal{T}_{c}\left[\spded\right](\zf, t) 
        \!= \spded(\zf,t) - \textstyle\int_0^{\zf} \bs{K}(\zf,\zeta,t)\spded(\zeta,t) \dop{\zeta}
\end{equation}
with the kernel $\bs{K}(\zf,\zeta,t)\in\mathbb{R}^{2\times2}$ is used to map
system \eqref{eq:sys_1d_fix_dec}
into the target system
\begin{subequations}
    \label{eq:target_sys}
    \begin{fleqn}
        \begin{align}
            \label{eq:target_sys_ode}
            \dot\sode(t) &\!=\! \bar{\bs{F}} \sode(t) 
                \!+\! \bs{S}(t) \partial_{\zf}\spdet(0,t)
                \\
            \label{eq:target_sys_pde}
            \begin{split}
                \partial_t \spdet(\zf,t) &\!=\!
                    \bs{\Lambda}(t) \partial_{\zf}^2 \spdet(\zf,t)
                    \!-\! \bs{M} \spdet(\zf,t)
        			\!-\! \bs{D}(\zf, t)\partial_{\zf} \spdet(0, t)
        			\raisetag{0.92\baselineskip}
            \end{split}
            \\	
            \label{eq:target_sys_bc0}
            \spdet(0, t) &\!= \bs{0} \\
            \label{eq:target_sys_bc1}
    		\partial_{\zf}\spdet(1, t) &\!= \bs{0},
        \end{align}
    \end{fleqn}
\end{subequations}
with the matrix $\bs{M} = \diag\left(\mu_{1}, \mu_{2}\right)$ containing the design parameters.
For that, as can easily be verified by inserting \eqref{eq:backstepp_trafo} into \eqref{eq:target_sys} and comparing to \eqref{eq:sys_1d_fix_dec}, $\bs{K}(\zf,\zeta,t)$ has to satisfy the kernel equations
\begin{subequations}
    \label{eq:kernel_equations}
    \begin{fleqn}
        \begin{align}
            \label{eq:k_pde}
            \partial_t \bs{K}(\zf,\zeta,t) &= 
                \bs{\Lambda}(t)\partial_{\zf}^2\bs{K}(\zf,\zeta,t)
                -\partial_{\zeta}^2\bs{K}(\zf,\zeta,t)\bs{\Lambda}(t)
                \nonumber\\&\;\quad
                - \bs{M}\bs{K}(\zf,\zeta,t)
                - \bs{K}(\zf,\zeta,t)\bs{A}(\zeta,t) 
            \raisetag{0.92\baselineskip}
            \\
            \label{eq:k_z0}
            \bs{K}(\zf, 0, t)\bs{\Lambda}(t) &= -\mathcal{T}_{c}[\bs{\bar C}](\zf,t)
                - \bs{D}(\zf,t)
            \\
            \begin{split}
                -(\bs{A}(\zf, t) \!+\! \bs{M}) &= \bs{\Lambda}(t)\left(
                    \tdiff{\zf}\bs{K}(\zf,\zf,t) 
                    + \partial_{\zf}\bs{K}(\zf,\zf,t) 
                    \right)
                    \\&\;\quad
                    + \partial_{\zeta}\bs{K}(\zf,\zf,t)\bs{\Lambda}(t)
                    \raisetag{2em}
            \end{split}
            \\
            \label{eq:k_zz}
            \bs{\Lambda}(t)\bs{K}(\zf,\zf,t) &= \bs{K}(\zf,\zf,t)\bs{\Lambda}(t),
        \end{align}
    \end{fleqn}
\end{subequations}
with the input given by
\begin{equation}
    \begin{split}
        \label{eq:control_law}
        \inpev &= \bs{Q}^{-1}(t) \Big[
          \textstyle\int_0^1 \partial_{\zf}\bs{K}(1,\zeta,t)\spded(\zeta,t) \dop{\zeta}
            \\&\;\;\;\quad
            + \left(\bs{K}(1,1,t) - \bs{B}(t)\right) \spded(1,t)
            - \bar{\bs{P}}(t)\sode(t)
        \Big].
    \end{split}
\end{equation}

To the best of the authors' knowledge, the kernel equations \eqref{eq:kernel_equations} have not yet been studied in the literature.
However, the solution of similar kernel equations is discussed in \cite{Kerschbaum2020,Kerschbaum2021}.
Therein, it is shown that \eqref{eq:kernel_equations} admits a piecewise classical solution of some Gevrey order in the case $\dot{\bs{\Lambda}}(t)=\bs{0}$ and $\bar{\bs{C}}(\zf,t)=\bs{0}$.
Obviously, a thorough discussion of the solution for $\dot{\bs{\Lambda}}(t)\neq\bs{0}$ and $\bar{\bs{C}}(\zf,t)\neq\bs{0}$ is beyond the scope of this contribution.
Still, from the approach in \cite{Kerschbaum2020,Kerschbaum2021} that uses a transformation of the kernel equations into integral form, fixed-point iterations and the method of successive iteration, it is clear that $\dot{\bs{\Lambda}}(t)\neq\bs{0}$ only increases the notation complexity,
even more so in view of Remark \ref{rem:diffsort}.
The case $\bar{\bs{C}}(z,t)\neq\bs{0}$ is more challenging.
In fact, \eqref{eq:k_z0} is only a \gls{bc} for those elements $k_{i,j}(\zf,\zeta,t)$ of $\bs{K}(\zf,\zeta,t)$ where $\lambda_i(t)>\lambda_j(t)$,
with the remaining equations defining the elements of $\bs{D}(\zf,t)$.
Hence, based on Remark \ref{rem:diffsort}, different cases have to be considered which is omitted for brevity.


\subsection{Final result}

By the tracking controller \eqref{eq:control_law} with the dynamics \eqref{eq:dyn_ext}, the closed loop \eqref{eq:target_sys} is obtained.
As $\bs{D}(\zf, t)$ in \eqref{eq:target_sys_pde} is either a strictly lower or upper triangular matrix for any $t$, choosing $\mu_{\gen} > 0$, $i=1,2$ results in a cascade of two exponentially stable parabolic \glspl{pde}.
Thus, the Hurwitz property of $\bar{\bs{F}}$ guarantees exponential stability of the closed loop \eqref{eq:target_sys}.
By the invertibility of the transformations \eqref{eq:decoupl_trafo} and \eqref{eq:backstepp_trafo}, the states $\sode(t)$ and $\spde(\zf,t)$ converge to zero, too.
This implies tracking of \eqref{eq:reference} for values of $\ts(z,t)$, $\tl(z,t)$ and $\zi$ close to the reference, in view of the linearisation in Section \ref{sec:prelim_trans}.

%% file: simulation.tex
\section{Numerical Results}
\label{sec:simulation}

To check the feedback design, a \SI{30}{\hour} excerpt from a
typical growth process for a \gls{gaas} single crystal is examined.

\begin{figure}[t]
    \centering
    \includegraphics[width=\linewidth]{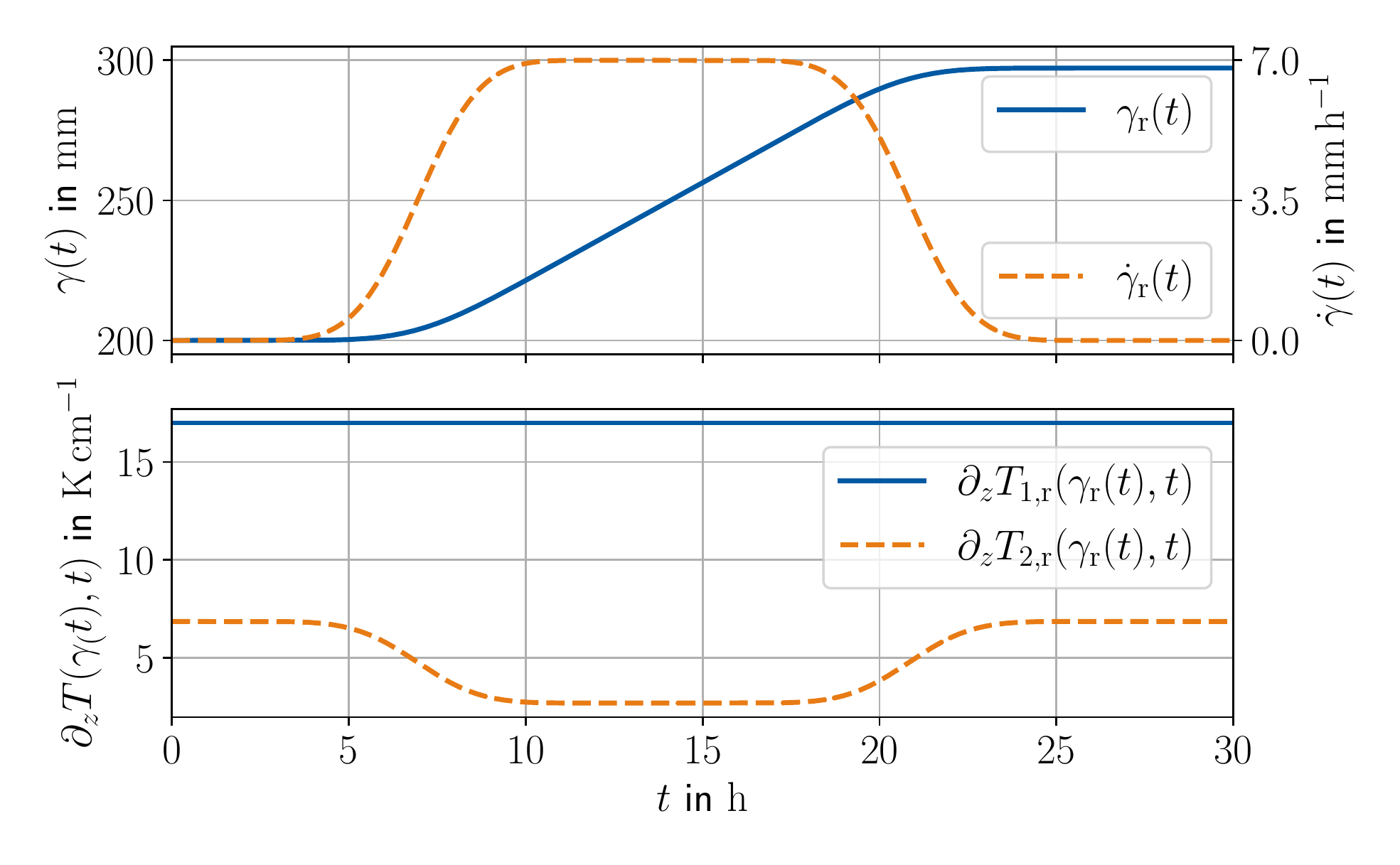}
    \caption{
       Solid lines represent the planned reference trajectories of the interface position $\zir$ and the temperature gradient {$\partial_z T_{1,\text{r}}(\zir,t)$} at the interface. Based on that, dashed lines indicate calculated trajectories for the growth rate $\dot\gamma_{\mathrm{r}}(t)$ of the crystal and the gradient {$\partial_z T_{2,\text{r}}(\zir,t)$} in the melt (implied by \eqref{eq:stefan_cond_sep}).
    }
    \label{fig:references}
\end{figure}

\subsection{Simulation Setup}
\label{sec:simsetup}

In accordance with the flatness-based feedforward design in Section~\ref{sec:feedforward} and \cite{Ecklebe2020}, reference trajectories are specified for the components of the flat output $\bs{y}(t)$, i.e.\ for the position $\zi$ of the interface and the temperature gradient $\partial_{\zf}T_{1}(\zi,t)$ thereat.
Specifically, an initial crystal size of
\SI{200}{\milli\metre} is assumed, from which a crystal of about \SI{300}{\milli\metre} is to be grown by smoothly transitioning
the growth velocity $\vi$ from initially \SI{0}{\milli\metre\per\hour}
to the nominal speed of \SI{7}{\milli\metre\per\hour} and back to zero.
Furthermore, to guarantee the desired crystal quality, 
a temperature gradient of \SI{17}{\kelvin\per\centi\meter} is to be held in the crystal at the interface at all times.
The reference trajectories corresponding to these specifications are given in Figure \ref{fig:references}, where a Gevrey order of $1.9$ is chosen for $\zir$.
Based on $\bs{y}_\text{r}(t)$, all other reference trajectories in \eqref{eq:reference} are calculated.

The simulation of the \gls{vgf} plant \eqref{eq:lin_heat_eqs} makes use of the numerical method in \cite[Sec.~3]{Ecklebe2021} with $64$ elements per phase.
The model parameters regarding the crucible and the material grown are taken from \cite[Tab.~I]{Ecklebe2022}.
Note that $\lambda_1(t)<\lambda_2(t)$ follows for these parameters and the chosen reference trajectories (cf.\ Remark~\ref{rem:diffsort}).
The initial conditions of \eqref{eq:lin_heat_eqs} are chosen such that initial errors of \SI{-10}{\milli\metre} for the position $\zi$ of the interface, $\SI{-3}{\milli\metre\per\hour}$ for the velocity $\vi$ and $\SI{5}{\kelvin\per\centi\metre}$ for the temperature gradient in the crystal are introduced (w.r.t.\ the specified reference trajectory in Figure \ref{fig:references}).

\subsection{Solution of Decoupling and Kernel Equations}

The implementation of the tracking controller \eqref{eq:control_law} (with the dynamics \eqref{eq:dyn_ext}) requires the solutions $\bs{N}(\zf,t)$ and $\bs{K}(\zf,\zeta,t)$ of \eqref{eq:decoup_equations} and \eqref{eq:kernel_equations}, respectively.
An approximated, numerical solution of the decoupling equations \eqref{eq:decoup_equations} is found by using a power series ansatz $\bs{N}(\zf,t)=\sum_{j=0}^{20}\bs{L}_j(t)\frac{\zf^j}{j!}$ of order $20$ and by choosing $\bs{K}(t)$ in \eqref{eq:def_barF} such that the \gls{ode} dynamics \eqref{eq:target_sys_ode} of the target system are governed by the Hurwitz matrix $\bar{\bs{F}} = \diag(\SI{2e-4}{\per\second},  \SI{2e-4}{\per\second})$.
The kernel equations \eqref{eq:kernel_equations} are solved in a simplified version that neglects the time derivative in \eqref{eq:k_pde}.
Thus, $\bs{K}(\zf,\zeta,t)$ is assumed to only slowly vary in time, which was confirmed through simulations.
Based on the chosen control gains $\mu_i = \SI{3e-4}{\per\second}$, $i=1,2$, the \texttt{ParabolicBackstepping}-toolbox (see \cite{KerschbaumCode}) solves the kernel equations \eqref{eq:kernel_equations} at $10$ equidistant points in time on a spatial grid of $101$ points.
Note that interpolating between
these $10$ quasistatic solutions (instead of solving the underlying time-varying problem \eqref{eq:kernel_equations}) induces an error in the evaluation of \eqref{eq:control_law}.
Both solutions $\bs{N}(\zf,t)$ and $\bs{K}(\zf,\zeta,t)$ are calculated offline.

\subsection{Results}
\label{sec:simresults}

\begin{figure}[t]
    \centering
    \includegraphics[width=\linewidth]{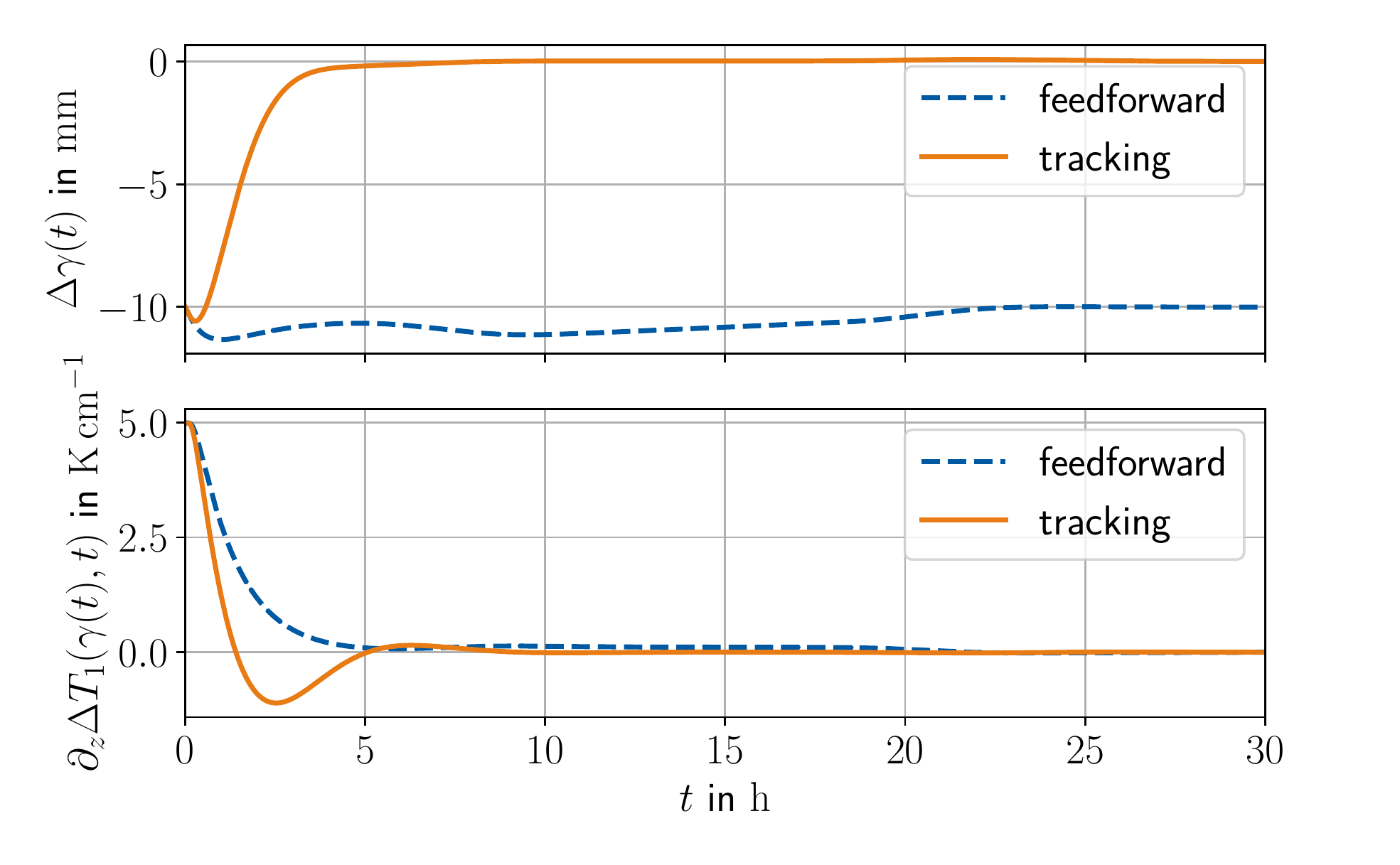}
    \caption{Errors $\dzi=\zi-\zir$ and $\partial_{\zf}\Delta T_1(\zi,t)=\partial_{\zf}T_1(\zi,t)-\partial_{\zf}T_{1,\text{r}}(\zir,t)$ of the flat output components resulting from the feedforward controller (dashed) and the tracking controller (solid).}
    \label{fig:ode_sim}
\end{figure}

\begin{figure}[t]
    \centering
    \includegraphics[width=\linewidth]{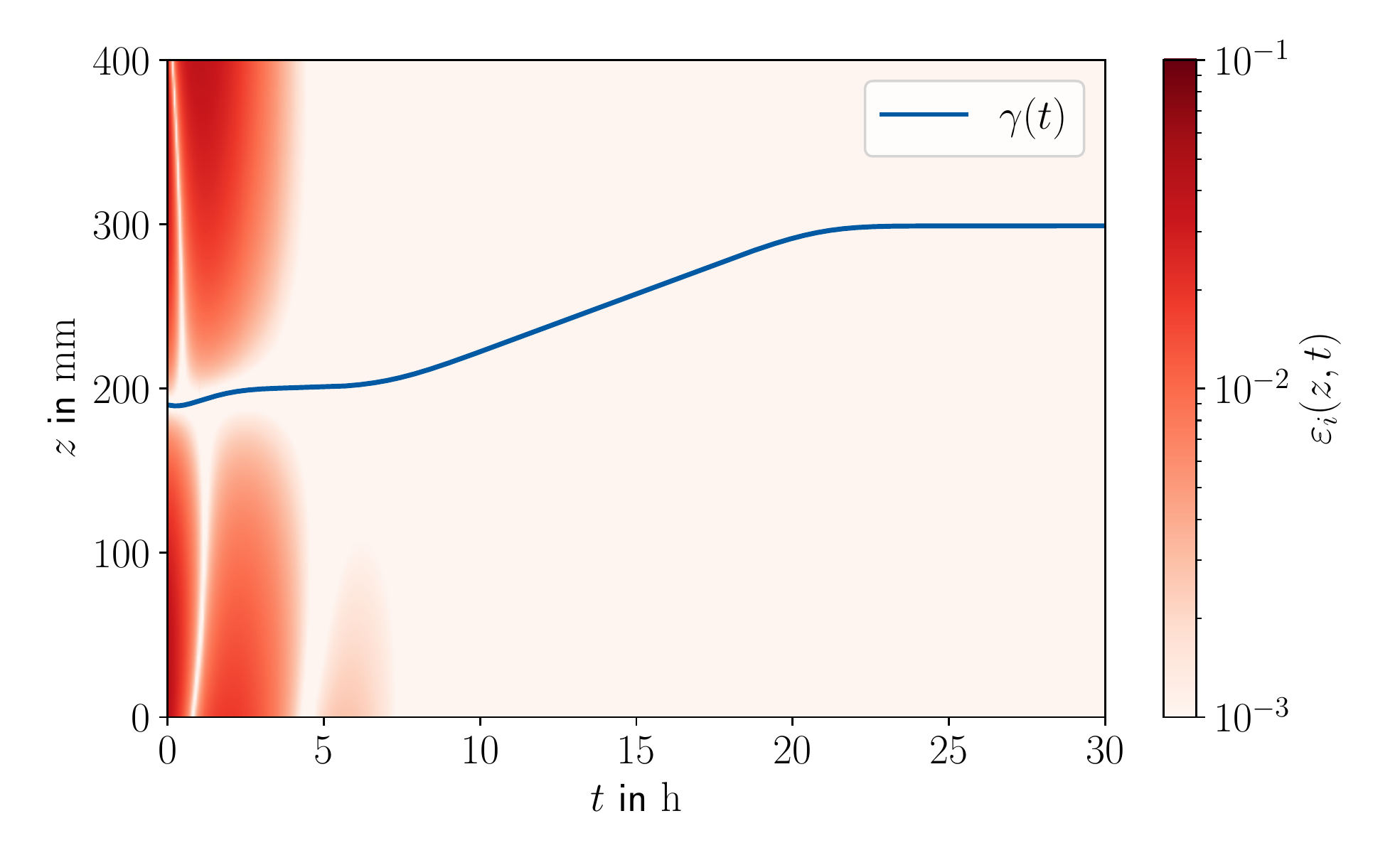}
    \caption{
        Logarithmic error $\varepsilon_i(\zf,t)\coloneqq\log\left|\frac{T_i(z,t)-T_{i,\text{r}}(z,t)}{T_{i,\text{r}}(z,t)}\right|$ of the joint temperature profiles $T_i(z,t)$ of crystal ($i=1$) and melt ($i=2$), relative to their references $T_{i,\text{r}}(z,t)$.
        Herein, both phases are separated by $\zi$.
    }
    \label{fig:pde_err_2d}
\end{figure}

Applying the tracking controller \eqref{eq:control_law} with the dynamics \eqref{eq:dyn_ext} to the \gls{vgf} plant \eqref{eq:lin_heat_eqs} yields the results in Figures \ref{fig:ode_sim} and \ref{fig:pde_err_2d}.
In particular, Figure \ref{fig:ode_sim} confirms that the controller tracks the references (see Figure \ref{fig:references}) for the position of the interface and the temperature gradient at the interface.
The initial errors converge to zero (or values close to that) in about \SI{5}{\hour}, which is also confirmed by the logarithmic plot of the relative error of the distributed temperature profile in crystal and melt in Figure \ref{fig:pde_err_2d}.
To better assess the control performance, Figure \ref{fig:ode_sim} gives the errors $\Delta\zi$ and $\partial_{\zf}\Delta T_{1}(\zi,t)$ if only a feedforward controller $\inps=\inprs$, $\inpl=\inprl$ is used.
With that, the error in the position of the interface cannot be compensated.
The temperature gradient at the interface converges, yet no steady-state accuracy is achieved during the main growth phase between \SI{10}{\hour} and \SI{18}{\hour} (see also Figure~\ref{fig:references}).
In contrast to the feedforward controller, the tracking controller creates a significant undershoot in the gradient.
This effect is necessary as the temperature gradient in the crystal has to be lowered below its reference to increase the growth velocity and catch up to the reference interface position $\zir$. 

%